%% file: A_Fibonacci_theorem_for_Collatz_trajectories.tex
\documentclass[12pt]{article}

\usepackage[utf8]{inputenc}
\usepackage[T1]{fontenc}
\usepackage{amsmath, amssymb, amsthm}
\usepackage[
pdftitle={A Fibonacci theorem for Collatz trajectories via modular graph structure},
pdfauthor={Manuel-Alejandro Reyes Jiménez},
pdfsubject={Collatz conjecture},
pdfkeywords={Collatz, 3x+1, Fibonacci, modular graph, binary encoding}
]{hyperref}
\usepackage{geometry}
\usepackage{enumitem}
\usepackage{titlesec}
\usepackage{mathrsfs}
\usepackage{graphicx}
\usepackage{caption}
\usepackage{float}
\usepackage{tikz}
\usepackage[table]{xcolor}
\usepackage{booktabs}
\usepackage{colortbl}
\usepackage{array}
\usepackage[activate=false]{microtype}
\usepackage{cite}

\usetikzlibrary{arrows.meta, positioning, shapes, calc,
                decorations.pathreplacing}

\newtheorem{theorem}{Theorem}[section]
\newtheorem{lemma}[theorem]{Lemma}
\newtheorem{proposition}[theorem]{Proposition}
\newtheorem{corollary}[theorem]{Corollary}

\theoremstyle{definition}
\newtheorem{definition}[theorem]{Definition}
\newtheorem{example}[theorem]{Example}

\theoremstyle{remark}
\newtheorem{remark}[theorem]{Remark}

\titleformat{\section}[block]{\Large\bfseries}{\thesection.}{0.5em}{}
\titleformat{\subsection}[block]{\large\bfseries}{\thesubsection.}{0.5em}{}

\title{A Fibonacci theorem for Collatz trajectories\\
       via modular graph structure}

\author{Manuel-Alejandro Reyes Jim\'enez}

\date{}

\begin{document}

\maketitle
\thispagestyle{empty}

% --- Abstract ---
\input{abstract_paper1.tex}



\bigskip

% --- Introduction ---
\input{intro_paper1.tex}

% --- Section 2: Binary and modular encodings ---
\input{sec2_new.tex}

% --- Section 3: Graph-theoretic structure of G ---
\input{sec3_new.tex}

% --- Section 4: Block decomposition and Fibonacci theorem ---
\input{sec4_block_structure.tex}

% --- Conclusions ---
\input{conclusions_paper1.tex}

% --- Acknowledgements ---
\section*{Acknowledgements}
The author thanks Jaume Franch Bullich (UPC, Barcelona),
PhD supervisor, for guidance and helpful discussions
throughout this work.

% --- Bibliography ---
\bibliographystyle{plain}
\bibliography{bibliografia_paper1}

\bigskip
\noindent
Manuel-Alejandro Reyes Jim\'enez\\
Departament de Matem\`atiques\\
Universitat Polit\`ecnica de Catalunya (UPC)\\
Barcelona, Spain\\
E-mail: \texttt{alejandro.reyes@upc.edu}\\
ORCID: \href{https://orcid.org/0009-0002-7216-9385}{0009-0002-7216-9385}

\end{document}

%% file: abstract_paper1.tex
\begin{abstract}
Let $T(n)=n/2$ if $n$ is even and $T(n)=(3n+1)/2$ if $n$ is odd.
We prove that for each $m\ge1$, exactly $F(m+1)$ odd integers~$n$ in
$\{1,\ldots,2^m\}$ have the property that none of the iterates
$T(n),T^2(n),\ldots,T^{m-1}(n)$ lies in the residue class $4\pmod6$,
where $F(m+1)$ is the $(m+1)$-th Fibonacci number; the proportion
decays at rate $(\varphi/2)^m$, $\varphi=(1+\sqrt{5})/2$.
Equivalently, these are the odd $n\le2^m$ for which no two consecutive
terms of $n,T(n),\ldots,T^{m-1}(n)$ are even.

The proof uses the directed graph $G$ of Collatz transitions modulo $6$
and its unique absorbing strongly connected component $G'=G[\{1,2,4,5\}]$.
Removing vertex $4$ from $G'$ yields a subgraph of spectral
radius $\varphi$, against $\rho(G')=2$; the Fibonacci count follows
from this spectral gap. We construct an explicit bijection
$\Psi_m:\{1,\ldots,6\cdot2^m\}\to\mathcal{P}_m(G)$ onto the
directed paths of length $m$ in $G$.

We further show that no vertex of $G'$ is dispensable: removing any
single vertex reduces the spectral radius strictly below $2$, with
hierarchy $1<\sqrt{2}<\varphi<2$. In particular, every positive cycle
of $T$ must visit residue class $2\pmod6$, and a flow
conservation identity forces this class to account for more
than $18\%$ of the steps in any such cycle.
\end{abstract}

%% file: intro_paper1.tex
\section{Introduction}
\label{sec:introduction}

The accelerated Collatz map
\[
T(n)=
\begin{cases}
n/2, & n\equiv0\pmod{2},\\[2pt]
(3n+1)/2, & n\equiv1\pmod{2},
\end{cases}
\]
defined over $\mathbb{Z}$, is conjectured to drive every positive integer
to the cycle $1\leftrightarrow2$, but no proof is known
\cite{Lagarias1985, Lagarias2006}. Probabilistic methods
\cite{Terras1976, Wirsching1998, Lagarias2006} describe average
trajectory behaviour; the deepest partial result is due to
Tao~\cite{Tao2022}, who showed that almost all positive integers
eventually reach arbitrarily small values.

The main result of this paper is the following. For each $m\ge1$,
among the odd integers~$n$ in $\{1,\ldots,2^m\}$, exactly $F(m+1)$
have the property that none of the iterates
$T(n),T^2(n),\ldots,T^{m-1}(n)$ lies in the residue
class~$4\pmod6$ --- equivalently, that no two consecutive terms of
$n,T(n),\ldots,T^{m-1}(n)$ are even
(Remark~\ref{rem:palabras-paridad}) --- where $F(m+1)$ is the
$(m+1)$-th Fibonacci number (Theorem~\ref{thm:fibonacci}); the
proportion decays at the exponential rate $(\varphi/2)^m$,
$\varphi=(1+\sqrt{5})/2$. The connection is structural: the entries
of the transfer matrix of~$\mathcal{H}_4$ are themselves Fibonacci
numbers (Proposition~\ref{prop:fibonacci-matrix}).

Fibonacci numbers have appeared previously in connection with the Collatz map.
Albert, Gudmundsson, and Ulfarsson~\cite{AlbertGudmundssonUlfarsson2015}
showed that Collatz permutations of length up to~$14$ are enumerated by the
Fibonacci numbers, a result that arises from counting word types over $\{u,d\}$.
Winkler~\cite{Winkler2014} conjectured, and verified computationally for
$2\le h,t\le50$, that the admissible residue classes of certain finite
subsequences of a Collatz trajectory are counted by the Fibonacci numbers;
this was proved in~\cite{Winkler2026}, where the count is obtained from
binary parity words avoiding the factor~$11$.

The object counted here is different: integers, classified by the avoidance
of a single residue class modulo~$6$, rather than residue classes determined
by a local cut condition on the parity word. The counting mechanism, however,
is of the same finite-state type, as observed in~\cite{Winkler2026}: a local
restriction on the orbit is encoded by paths in a finite state space and then
evaluated by a transfer matrix. What distinguishes the present approach is
that its states are residues modulo~$6$ rather than local parity data. The
golden ratio $\varphi$ accordingly appears as the Perron--Frobenius
eigenvalue of an explicit subgraph~$\mathcal{H}_4$ of the modular
graph~$G'$, and the same graph carries the spectral hierarchy and the cycle
constraints described below, which have no counterpart at the level of parity
words alone.

The proof rests on two parallel objects, one binary and one modular.
On the binary side, we use the parity encoding
$C^m(n)=(c_1,\ldots,c_m)$ of the first $m$ iterates and its associated
\emph{corrective residue} $R(C^m(n))$ (Lemma~\ref{lem:residuo-estructural}).
The key identity is $T^m(n)=(3^s n + R)/2^m$. Since $C^m(n)$
depends on~$n$ only through its residue modulo~$2^m$, $R$ is constant
on residue classes modulo~$2^m$. Algebraically
equivalent identities appear in Terras~\cite{Terras1976},
B\"ohm--Sontacchi~\cite{BohmSontacchi1978}, and
Lagarias~\cite{Lagarias1985}.
The map $\Phi_m:n\mapsto C^m(n)$
is a bijection from $\{1,\ldots,2^m\}$ onto $\{0,1\}^m$, due to
Everett~\cite{Everett1977} and
Bernstein--Lagarias~\cite{BernsteinLagarias1996}, and satisfies
periodicity $C^m(n+2^m)=C^m(n)$ and antisymmetry
$c_{m+1}(n+2^m)=\overline{c_{m+1}(n)}$. Together, these properties
allow the entire analysis to be conducted at the level of binary
encodings, without reference to the specific value of~$n$.

The modular analogue encodes the dynamics modulo~$6$ --- the smallest
modulus that simultaneously distinguishes parity and the image of the
odd branch of~$T$ --- via a directed graph $G$ with vertex set
$\{0,1,2,3,4,5\}$ and twelve edges. Modular Collatz graphs for general
modulus~$N$ are defined in~\cite{KarrasWeger2026}, where the case $N=6$
appears as an explicit example. The present paper develops a different
treatment: we construct a bijection
$\Psi_m:\{1,\ldots,6\cdot2^m\}\to\mathcal{P}_m(G)$ onto directed paths
of length~$m$ in~$G$, satisfying periodicity and antisymmetry properties
that parallel those of the binary encoding~$\Phi_m$. To our knowledge,
this bijection and its structural properties do not appear in prior work.

The graph $G$ has a unique absorbing strongly connected component
$G'=G[\{1,2,4,5\}]$ (Proposition~\ref{prop:scc}); odd integer orbits
are confined to~$G'$ from level~$2$ onwards
(Corollary~\ref{cor:subgrafo-efectivo}). The Fibonacci
count arises from a spectral gap: the avoiding subgraph
$\mathcal{H}_4=G'[\{1,2,5\}]$, obtained by removing vertex~$4$, has
$\rho(\mathcal{H}_4)=\varphi$, while $\rho(G')=2$, and the ratio
$\varphi/2<1$ is the exact exponential decay rate (Corollary~\ref{cor:fibonacci-densidad}). A modular path
theorem (Theorem~\ref{thm:camino-modular}) describes the action of
the displacement~$+2^m$ on modular paths, with the
consequence that odd cycles of length~$m$ are separated by at least
$6\cdot2^m$ in any arithmetic progression of common difference~$2^m$
(Corollary~\ref{cor:ciclos}).

More generally, removing any vertex $v$ from~$G'$ reduces the spectral
radius strictly below~$2$. The exact values,
$\rho(\mathcal{H}_2)=1$, $\rho(\mathcal{H}_5)=\sqrt{2}$,
$\rho(\mathcal{H}_4)=\rho(\mathcal{H}_1)=\varphi$,
form a hierarchy $1<\sqrt{2}<\varphi<2$ governing the exponential
decay of the corresponding avoiding densities
(Theorem~\ref{thm:indispensabilidad}). In particular, no positive
cycle of~$T$ can avoid vertex~$2$ (Corollary~\ref{cor:visits-v2}). A flow conservation identity
in~$G'$ further forces vertex~$2$ to account for more than
$m(1-\log_32)/2\approx0.1845\,m$ steps in any positive cycle of
length~$m$ (Corollary~\ref{cor:quantitative}).

The paper is organized as follows: Section~\ref{sec:binary} establishes the binary encoding~$\Phi_m$ and
the modular bijection~$\Psi_m$, together with their shared periodicity
and antisymmetry. Section~\ref{sec:graph} determines the graph-theoretic
and spectral structure of~$G$. The main theorems are proven in
Section~\ref{sec:paths}, using the results of
Sections~\ref{sec:binary} and~\ref{sec:graph}.
Section~\ref{sec:conclusions} states open questions
and discusses connections to ongoing work.

%% file: sec2_new.tex
\section{Binary and modular encodings}
\label{sec:binary}

We establish the two parallel bijections $\Phi_m$ and $\Psi_m$,
together with their shared periodicity and antisymmetry properties,
on which the proofs of Section~\ref{sec:paths} depend.

To each integer $n$ and level $m\ge1$ we associate the \emph{binary encoding}
\[
C^m(n)=(c_1,\ldots,c_m)\in\{0,1\}^m,\qquad c_i=T^{i-1}(n)\bmod2,
\]
recording the parities of $n,T(n),\ldots,T^{m-1}(n)$. The \emph{Hamming weight}
$s=s(C^m(n))=\sum_{i=1}^m c_i$ counts the odd steps, and
$r(i)=\sum_{j=1}^i c_j$ denotes the partial weight up to position~$i$.

We define the \emph{corrective residue} by
\begin{equation}\label{eq:residuo-def}
R(C^m(n)):=2^m\,T^m(n)-3^s\,n,
\end{equation}
which satisfies the structural identity
\begin{equation}\label{eq:iteracion}
T^m(n)=\frac{3^s\,n+R(C^m(n))}{2^m}.
\end{equation}
Algebraically equivalent expressions appear in \cite{Terras1976,
BohmSontacchi1978, Lagarias1985}; in particular,
\cite[Prop.~4]{BohmSontacchi1978} gives, for the accelerated map
$u(x)=(3x+1)/2$, a closed form whose numerator coincides with
$R(C^m(n))$ under the standard correspondence between gap-sequences
and binary encodings. Since $C^m(n)$ depends on~$n$ only through its
residue modulo~$2^m$ (Proposition~\ref{prop:biy-fractal}(B,P) below),
$R$ is constant on residue classes modulo~$2^m$; the bijection
$\Phi_m$ then allows the analysis to be conducted entirely at the
level of encodings, with no reference to the specific value of~$n$,
which drives the parallel modular construction below.

\begin{lemma}%
\label{lem:residuo-estructural}
For each $m\ge1$, the corrective residue depends solely on $C^m(n)=(c_1,\ldots,c_m)$ and is
given explicitly by
\begin{equation}\label{eq:residuo-correctivo}
R(C^m(n))=\sum_{i=1}^m 3^{\,s-r(i)}\,2^{\,i-1}\,c_i.
\end{equation}
Moreover, it satisfies the recursive law:
\[
R(C^{m+1}(n))=
\begin{cases}
R(C^m(n)), & c_{m+1}=0,\\
3\,R(C^m(n))+2^m, & c_{m+1}=1.
\end{cases}
\]
\end{lemma}

\begin{proof}
By induction on~$m$. The base case $m=1$ is immediate: $c_1=0$ gives
$R=0$, and $c_1=1$ gives $R=1$, both matching the formula.

For the inductive step, suppose
$R(C^m(n))=\sum_{i=1}^m 3^{\,s-r(i)}\,2^{\,i-1}\,c_i$
holds at level~$m$. We prove it at level~$m+1$.

\emph{Case $c_{m+1}=0$.} Then $s'=s$, $r'(i)=r(i)$ for all $i\le m+1$,
and $c_{m+1}=0$, so the new term in the sum vanishes:
\[
\sum_{i=1}^{m+1}3^{\,s'-r'(i)}\,2^{\,i-1}\,c_i
=\sum_{i=1}^{m}3^{\,s-r(i)}\,2^{\,i-1}\,c_i=R(C^m).
\]
On the other hand, $T^{m+1}(n)=T^m(n)/2$, so
$R(C^{m+1})=2^{m+1}T^{m+1}(n)-3^{s'}n=2^mT^m(n)-3^sn=R(C^m)$.
Both sides agree.

\emph{Case $c_{m+1}=1$.} Then $s'=s+1$ and $r'(i)=r(i)$ for $i\le m$,
$r'(m+1)=s+1$. The new term at position $m+1$ contributes
$3^{s'-r'(m+1)}\,2^m\cdot1=3^0\cdot2^m=2^m$.
Each existing term at position $i\le m$ acquires an extra factor of~$3$,
since $s'-r'(i)=(s+1)-r(i)=(s-r(i))+1$. Therefore:
\begin{align*}
\sum_{i=1}^{m+1}3^{\,s'-r'(i)}\,2^{\,i-1}\,c_i
&=\sum_{i=1}^{m}3^{\,(s-r(i))+1}\,2^{\,i-1}\,c_i+2^m=3\sum_{i=1}^{m}3^{\,s-r(i)}\,2^{\,i-1}\,c_i+2^m\\
&=3R(C^m)+2^m.
\end{align*}
On the other hand, $T^{m+1}(n)=(3T^m(n)+1)/2$, so
\begin{align*}
R(C^{m+1})&=2^{m+1}T^{m+1}(n)-3^{s+1}n=2^m(3T^m(n)+1)-3\cdot3^sn\\
&=3(2^mT^m(n)-3^sn)+2^m=3R(C^m)+2^m.
\end{align*}
Both sides equal $3R(C^m)+2^m$, completing the induction. The recursive
law is the $c_{m+1}$-case split just established.
\end{proof}

\begin{example}
For $n=11$, $m=4$: iterates $11\to17\to26\to13\to20$,
encoding $C^4(11)=(1,1,0,1)$, $s=3$. The formula gives
$R=9\cdot1+6\cdot1+0+1\cdot8=23$, and indeed
$T^4(11)=(27\cdot11+23)/16=320/16=20$.
\end{example}

Equation~\eqref{eq:iteracion} implies three exact properties of the encoding
used throughout Section~\ref{sec:paths}. Bijectivity~(B) and
periodicity~(P) are due to Everett~\cite{Everett1977} and
Bernstein--Lagarias~\cite{BernsteinLagarias1996};
antisymmetry~(A) does not appear to have been stated explicitly in prior work.

\begin{proposition}[Bijectivity~(B) and Periodicity~(P)]\label{prop:biy-fractal}
For each $m\ge1$:
\begin{itemize}
  \item[\emph{(B)}] The map $n\mapsto C^m(n)$ is a bijection from
    $\{1,\ldots,2^m\}$ onto $\{0,1\}^m$. Equivalently,
    $C^m(n_1)=C^m(n_2)\Leftrightarrow n_1\equiv n_2\pmod{2^m}$.
  \item[\emph{(P)}] The encoding has exact period~$2^m$:
    $C^m(n+2^m)=C^m(n)$ for all $n\in\mathbb{Z}$.
\end{itemize}
\end{proposition}

\begin{proof}
We prove~(P) first, as~(B) depends on it.

\emph{Periodicity~(P).} By induction on~$m$.

For $m=1$: $c_1(n)=n\bmod2$ has period~$2=2^1$, so $c_1(n+2)=c_1(n)$.

Assume (P) holds at level~$m$, i.e.\ $C^m(n+2^m)=C^m(n)$ for all~$n$.
We prove it at level~$m+1$. Since $2^{m+1}$ is a multiple of~$2^m$,
the inductive hypothesis gives $C^m(n+2^{m+1})=C^m(n)$, so
$c_i(n+2^{m+1})=c_i(n)$ for $i=1,\ldots,m$. It remains to show
$c_{m+1}(n+2^{m+1})=c_{m+1}(n)$. By~\eqref{eq:iteracion} applied at
level~$m$,
\[
T^m(n+2^{m+1})=\frac{3^s(n+2^{m+1})+R}{2^m}=T^m(n)+2\cdot3^s,
\]
where $s$ and $R$ are the same for $n$ and $n+2^{m+1}$ by the inductive
hypothesis. Since $2\cdot3^s$ is even, $T^m(n+2^{m+1})$ and $T^m(n)$ have
the same parity, hence $c_{m+1}(n+2^{m+1})=c_{m+1}(n)$.

\emph{Bijectivity~(B).} Injectivity: if $C^m(n_1)=C^m(n_2)$, then $s$ and
$R$ are the same for both, so $T^m(n_1)=(3^sn_1+R)/2^m$ and
$T^m(n_2)=(3^sn_2+R)/2^m$. Subtracting gives $3^s(n_1-n_2)\equiv0\pmod{2^m}$.
Since $\gcd(3^s,2^m)=1$, we have $n_1\equiv n_2\pmod{2^m}$. As $n_1,n_2\in
\{1,\ldots,2^m\}$, this gives $n_1=n_2$. Surjectivity follows from
injectivity, since both domain and codomain have exactly $2^m$ elements.
\end{proof}

\noindent
We write $\Phi_m:\{1,\ldots,2^m\}\to\{0,1\}^m$ for this bijection
henceforth.

\begin{lemma}[Antisymmetry~(A)]\label{lem:antisimet}
For each $m\ge1$ and all $n\in\mathbb{Z}$,
\[
c_{m+1}(n+2^m)=\overline{c_{m+1}(n)},
\]
where $\overline{c}:=1-c$ denotes the binary complement of $c\in\{0,1\}$.
The shift $+2^m$ preserves the entire level-$m$ encoding and inverts
precisely the $(m+1)$-th entry.
\end{lemma}

\begin{proof}
By~(P), $C^m(n+2^m)=C^m(n)$, so both share the same weight~$s$
and residue~$R$ at level~$m$. Then \eqref{eq:iteracion} gives
\[
T^m(n+2^m)=\frac{3^s(n+2^m)+R}{2^m}=\frac{3^s n+R}{2^m}+3^s=T^m(n)+3^s.
\]
Since $3^s$ is odd, $T^m(n+2^m)$ and $T^m(n)$ have opposite parities;
hence $c_{m+1}(n+2^m)=\overline{c_{m+1}(n)}$.
\end{proof}

\begin{lemma}\label{lem:traslacion}
For all $n\in\mathbb{Z}$, $t\in\mathbb{Z}$, $k\ge0$,
\[
T^m(n\pm t\cdot2^{m+k})=T^m(n)\pm t\cdot2^k\cdot3^s.
\]
This generalises~(P) and is used in Section~\ref{sec:paths} to
control how $T^m$ responds to arithmetic displacements at scales
above~$2^m$.
\end{lemma}

\begin{proof}
By~(P), $C^m(n\pm t\cdot2^{m+k})=C^m(n)$, so~\eqref{eq:iteracion}
gives $T^m(n\pm t\cdot2^{m+k})=(3^s(n\pm t\cdot2^{m+k})+R)/2^m
=T^m(n)\pm t\cdot2^k\cdot3^s$.
\end{proof}

A \emph{cycle of length~$m$} of~$T$ is an integer $n$ such that
$T^m(n)=n$ and $T^j(n)\ne n$ for $1\le j<m$. A \emph{positive cycle}
is a cycle with $n>0$.

\begin{lemma}\label{lem:cycle-eq}
For $m\ge1$ and $s\ge1$, the condition $T^m(n)=n$ is equivalent to
\begin{equation}\label{eq:ciclos}
n=\frac{R(C^m(n))}{2^m-3^s}.
\end{equation}
Note that $2^m\neq3^s$ whenever $s\ge1$, since $2$ and $3$ are
multiplicatively independent.
\end{lemma}

\begin{proof}
Equation~\eqref{eq:iteracion} gives $T^m(n)=n$ if and only if
$3^sn+R=2^mn$, i.e.\ $n(2^m-3^s)=R$. Since $s\ge1$ implies
$2^m\neq3^s$, we may divide to get $n=R/(2^m-3^s)$.
\end{proof}

\begin{remark}
The algebraic content of~\eqref{eq:ciclos} is due to
B\"ohm--Sontacchi~\cite{BohmSontacchi1978}; the form here follows
from Lemma~\ref{lem:residuo-estructural} by setting $T^m(n)=n$
in~\eqref{eq:iteracion}.
The hypothesis $s\ge1$ is necessary: if $s=0$, every iterate
$T(n),\ldots,T^{m-1}(n)$ is even, so $T^j(n)=n/2^j$ for
$j=1,\ldots,m$, which requires $2^m\mid n$. Then $T^m(n)=n/2^m<n$
for $n>0$; in particular, no positive integer satisfies $T^m(n)=n$
with $s=0$.
\end{remark}

The modular analogue encodes the dynamics modulo~$6$ via the graph $G$,
whose transitions follow from the general modular Collatz framework
of~\cite{KarrasWeger2026}. We construct
a bijection $\Psi_m:\{1,\ldots,6\cdot2^m\}\to\mathcal{P}_m(G)$,
where $\mathcal{P}_m(G)$ denotes the set of directed paths of length~$m$
in~$G$ (Definition~\ref{def:camino-modular}), satisfying the same
periodicity (Lemma~\ref{lem:periodicidad-modular})
and antisymmetry (see Remark below) as $\Phi_m$,
at scale $6\cdot2^m$.

The transitions of $T$ modulo~$6$ take exactly two values in each
residue class:

\begin{proposition}\label{prop:transiciones-mod6}
For each $a\in\{0,1,2,3,4,5\}$, the values of $T(n)\bmod6$
when $n\equiv a\pmod6$ are:
\[
\begin{aligned}
0&\to0,\;0\to3,\quad
1\to2,\;1\to5,\quad
2\to1,\;2\to4,\\
3&\to2,\;3\to5,\quad
4\to2,\;4\to5,\quad
5\to2,\;5\to5.
\end{aligned}
\]
Each residue class has exactly two possible images under~$T$.
\end{proposition}

\begin{proof}
Write $n=6k+a$. For even~$a$, $T(n)=(6k+a)/2=3k+a/2$, which takes
two values modulo~$6$ according to the parity of~$k$. For odd~$a$,
$T(n)=(3(6k+a)+1)/2=9k+(3a+1)/2\equiv 3k+(3a+1)/2\pmod6$, again
with two values. Evaluating each case:
\[
\begin{array}{c|c|c|c}
a & T(6k+a) & k\text{ even} & k\text{ odd} \\\hline
0 & 3k      & 0 & 3 \\
1 & 9k+2    & 2 & 5 \\
2 & 3k+1    & 1 & 4 \\
3 & 9k+5    & 5 & 2 \\
4 & 3k+2    & 2 & 5 \\
5 & 9k+8    & 2 & 5
\end{array}
\]
The twelve transitions follow.
\end{proof}

\begin{definition}[Modular Collatz graph]\label{def:grafo-modular}
The \emph{modular Collatz graph} is the directed graph $G=(V,E)$ with
vertex set $V=\{0,1,2,3,4,5\}$ and edge set
\begin{align*}
E=&\{(0,0),(0,3),(1,2),(1,5),(2,1),(2,4), (3,2),(3,5),(4,2),(4,5),\\
&(5,2),(5,5)\},
\end{align*}
the twelve transitions of Proposition~\ref{prop:transiciones-mod6}.
We write $a\to b$ for the edge $(a,b)\in E$.
The graph is depicted in Figure~\ref{fig:grafo-modular}.
\end{definition}

\begin{figure}[H]
\centering
\begin{tikzpicture}[
    >=Stealth,
    node distance=2.5cm,
    every node/.style={circle, draw, minimum size=10mm, font=\large},
    every path/.style={semithick},
    every loop/.style={looseness=6}
]
  \node (A) at (0,0) {0};
  \node (B) at (2,0) {3};
  \node (C) at (4,2) {5};
  \node (D) at (4,-2) {2};
  \node (E) at (5,0) {1};
  \node (F) at (7,0) {4};

  \draw[->] (A) edge[out=60, in=120, looseness=6] (A);
  \draw[->] (A) -- (B);
  \draw[->] (B) -- (C);
  \draw[->] (B) -- (D);
  \draw[->] (C) edge[out=60, in=120, looseness=6] (C);
  \draw[->] (C) -- (D);
  \draw[->, bend left=6] (D) to (E);
  \draw[->, bend left=6] (E) to (D);
  \draw[->] (E) -- (C);
  \draw[->] (F) -- (C);
  \draw[->, bend left=6] (D) to (F);
  \draw[->, bend left=6] (F) to (D);
\end{tikzpicture}
\caption{The modular Collatz graph $G$. Vertices~$0$ and~$3$ are transient;
the absorbing component is $G'=G[\{1,2,4,5\}]$.}
\label{fig:grafo-modular}
\end{figure}

\begin{definition}[Transition map]\label{def:tau}
The \emph{transition} performed by~$n$ is the edge
\[
\tau(n):=(n\bmod6,\;T(n)\bmod6)\in E.
\]
\end{definition}

The orbit of an integer $n$ under $T$ traces a walk in $G$, recorded
by its \emph{modular encoding}:

\begin{definition}[Modular encoding]\label{def:codificacion-mod-6}
For each integer~$n$ and level $m\ge1$, the \emph{modular encoding
of length~$m$} is
\[
H^m(n)=(h_1,\ldots,h_m),\qquad h_i\equiv T^{i-1}(n)\pmod{6},
\quad 1\le i\le m.
\]
The relation $c_i\equiv h_i\pmod{2}$ recovers the binary encoding.
\end{definition}

We formalise this walk as a path in $G$:

\begin{definition}[Modular path of length~$m$ and set $\mathcal{P}_m(G)$]%
\label{def:camino-modular}
Let $m\ge1$. A \emph{modular path of length~$m$} in~$G$ is a sequence of
$m$ consecutive edges
\[
\gamma=(e_1,e_2,\ldots,e_m)\in E^m,
\]
where $E$ is the edge set of Definition~\ref{def:grafo-modular}, and
each $e_i=(a_i,b_i)$ satisfies $b_i=a_{i+1}$ for all $1\le i\le m-1$
(the terminal vertex of~$e_i$ equals the initial vertex of~$e_{i+1}$). We denote by $\mathcal{P}_m(G)$ the set of all directed
paths of length~$m$ in~$G$.

The modular path \emph{induced by~$n$} at level~$m$ is
\[
\mathcal{G}_m(n):=\bigl(\tau(n),\;\tau(T(n)),\;\ldots,\;\tau(T^{m-1}(n))\bigr)
\in\mathcal{P}_m(G).
\]
Note that $\mathcal{G}_m(n)$ records the edge sequence
$(h_i,h_{i+1})$, $i=1,\ldots,m$, and is thus equivalent to the vertex
sequence $H^{m+1}(n)=(h_1,\ldots,h_{m+1})$.
\end{definition}

\begin{theorem}[Modular bijection]\label{thm:biyeccion-modular}
For every $m\ge1$ there is a bijection
\[
\Psi_m:\{1,\dots,6\cdot2^m\}\longrightarrow\mathcal{P}_m(G).
\]
In particular, $\mathcal{P}_m(G)$ has exactly $6\cdot2^m$ elements.
\end{theorem}

\begin{proof}
Call an edge \emph{of type~$\varepsilon\in\{0,1\}$} if its terminal
vertex has parity~$\varepsilon$, i.e.\ type~$0$ if the next iterate
is even and type~$1$ if it is odd. Since the vertices
of~$G$ are the residue classes $V=\{0,1,2,3,4,5\}$, each vertex
$a\in V$ has exactly two outgoing edges; by
Proposition~\ref{prop:transiciones-mod6} the two successors of~$a$
have opposite parities, so one outgoing edge has type~$0$ and the
other has type~$1$. Hence for each $a\in V$ and $\varepsilon\in\{0,1\}$
there is a unique edge of type~$\varepsilon$ leaving~$a$; we
denote it $e(a,\varepsilon)$.

\emph{Construction of~$\Psi_m$.}
Given $k\in\{1,\dots,6\cdot2^m\}$, perform the Euclidean division
\[
k-1=a\cdot2^m+r,\qquad a\in\{0,\dots,5\},\quad r\in\{0,\dots,2^m-1\},
\]
and write $r$ in base-$2$ representation with exactly $m$ digits:
$r=\sum_{i=1}^m\varepsilon_i\,2^{m-i}$, $\varepsilon_i\in\{0,1\}$.
Define $\Psi_m(k):=(e_1,\dots,e_m)\in\mathcal{P}_m(G)$ by
\[
e_1:=e(a,\varepsilon_1),\qquad
e_i:=e(b_{i-1},\varepsilon_i)\quad\text{for }i\ge2,
\]
where $b_{i-1}$ denotes the terminal vertex of~$e_{i-1}$.
This is well defined: each $e_i$ is determined by the terminal vertex of the
previous edge and the type $\varepsilon_i$, both uniquely specified.

\emph{Injectivity.}
Suppose $\Psi_m(k)=\Psi_m(k')$, so $(e_1,\dots,e_m)=(e_1',\dots,e_m')$.
Write $e_1=(a,b_1)$ and $e_1'=(a',b_1')$.
From $e_1=e_1'$ we read off $a=a'$.
Since $e(a,0)\ne e(a,1)$ (distinct types give distinct edges),
$e_1=e_1'$ also forces $\varepsilon_1=\varepsilon_1'$.
By induction: if $\varepsilon_j=\varepsilon_j'$ for all $j<i$, then
the terminal vertex of $e_{i-1}$ equals that of $e_{i-1}'$, and from
$e_i=e_i'$ the same argument gives $\varepsilon_i=\varepsilon_i'$.
Hence $\varepsilon_i=\varepsilon_i'$ for all~$i$, so $r=r'$ and $k=k'$.

\emph{Surjectivity.}
Let $(e_1,\dots,e_m)\in\mathcal{P}_m(G)$ be any path. Write $e_1=(a,b_1)$
and set $a\in\{0,\dots,5\}$ as its initial vertex. For each~$i$, write
$e_i=(a_i,b_i)$ and define $\varepsilon_i\in\{0,1\}$ as the unique type such that
$e_i=e(a_i,\varepsilon_i)$; this is well defined
since $e(a_i,0)$ and $e(a_i,1)$ are the only two edges leaving~$a_i$.
Set $r:=\sum_{i=1}^m\varepsilon_i\,2^{m-i}\in\{0,\dots,2^m-1\}$
and $k:=a\cdot2^m+r+1$. Then $k\in\{1,\dots,6\cdot2^m\}$
(since $a\le5$ and $r\le2^m-1$), and by construction
$\Psi_m(k)=(e_1,\dots,e_m)$.
\end{proof}

\begin{lemma}\label{lem:periodicidad-modular}
Let $m\ge1$. For all $n\in\mathbb{Z}$,
\[
\mathcal{G}_m\bigl(n+6\cdot2^m\bigr)=\mathcal{G}_m(n),
\]
where $\mathcal{G}_m(n)$ denotes the modular path induced by~$n$
(Definition~\ref{def:camino-modular}).
\end{lemma}

\begin{proof}
Set $n'=n+6\cdot2^m$. Since $6\cdot2^m$ is a multiple of~$6$, we have
$n'\equiv n\pmod6$, so $h_1(n')=h_1(n)$. Since $6\cdot2^m=3\cdot2^{m+1}$ is a
multiple of~$2^{m+1}$, Proposition~\ref{prop:biy-fractal}(P) applied at
level~$m+1$ gives $C^{m+1}(n')=C^{m+1}(n)$, so $c_i(n')=c_i(n)$ for
$i=1,\ldots,m+1$.

We now prove $h_i(n')=h_i(n)$ for $i=1,\ldots,m+1$ by induction on~$i$.
The base case $i=1$ holds by the above. For the inductive step, assume
$h_i(n')=h_i(n)$. By Proposition~\ref{prop:transiciones-mod6} the two
successors of any vertex have opposite parity, so $h_{i+1}$ is determined
by $h_i$ together with $c_{i+1}$ (the parity of~$T^i(n)$); since both agree
for $n$ and $n'$, we get $h_{i+1}(n')=h_{i+1}(n)$. As the $m$ edges
of~$\mathcal{G}_m$ are $(h_i,h_{i+1})$ for $i=1,\ldots,m$, this gives
$\mathcal{G}_m(n')=\mathcal{G}_m(n)$.
\end{proof}

\begin{remark}
The bijection~$\Psi_m$ satisfies the same antisymmetry as~$\Phi_m$,
at modular scale. Setting $\delta_m:=\operatorname{lcm}(6,2^m)=3\cdot2^m$ for $m\ge1$,
the displacement $+\delta_m$ preserves the first~$m$ levels of
$\mathcal{G}_{m+1}(n)$ and complements the edge at level~$m+1$:
\[
h_i(n+\delta_m)=h_i(n)\;\text{ for }1\le i\le m,
\qquad
h_{m+1}(n+\delta_m)=\overline{h_{m+1}(n)}.
\]
Here, for $h\in\{0,\ldots,5\}$, $\overline{h}:=h+3\pmod6$; since the two
successors of any vertex differ by $3\pmod6$, $\overline{h_{m+1}}$ is
precisely the other successor of~$h_m$, consistent with the binary
complement of Lemma~\ref{lem:antisimet} (complementing a vertex
complements its parity).
This follows from Lemma~\ref{lem:antisimet} and
Proposition~\ref{prop:transiciones-mod6}: the factor~$6$ in~$\delta_m$
preserves $h_1$, the factor~$2^m$ inverts $c_{m+1}$, and since each
vertex of~$G$ has exactly one even and one odd successor, the opposite
parities of $T^m(n)$ for $n$ and $n+\delta_m$ force the edge at
level~$m+1$ to be complemented.
\end{remark}

\begin{remark}[Structural parallel between $\Phi_m$ and $\Psi_m$]
The two bijections share the same three structural properties, at their
respective scales:
\begin{itemize}
  \item[\emph{(B)}] $\Phi_m$ bijects $\{1,\ldots,2^m\}$ onto $\{0,1\}^m$;
    $\Psi_m$ bijects $\{1,\ldots,6\cdot2^m\}$ onto $\mathcal{P}_m(G)$.
  \item[\emph{(P)}] Each has its own exact period, $2^m$ for $\Phi_m$ and
    $6\cdot2^m$ for $\Psi_m$:
    $C^m(n+2^m)=C^m(n)$ and $\mathcal{G}_m(n+6\cdot2^m)=\mathcal{G}_m(n)$.
  \item[\emph{(A)}] Both invert precisely the $(m+1)$-th entry under their
    respective antisymmetry displacements~$+2^m$ and~$+\delta_m$.
\end{itemize}
This parallel structure is not a coincidence: both bijections encode the
same dynamical information — one at the binary level, one at the modular
level — and the proofs of Sections~\ref{sec:graph} and~\ref{sec:paths}
exploit both simultaneously.
\end{remark}

%% file: sec3_new.tex
\section{Structure and spectrum of \texorpdfstring{$G$}{G}}
\label{sec:graph}

\begin{proposition}%
\label{prop:scc}
The graph~$G$ of Definition~\ref{def:grafo-modular} decomposes into
three strongly connected components:
\[
V_1=\{0\},\qquad V_2=\{3\},\qquad V_3=\{1,2,4,5\}.
\]
The components $V_1$ and $V_2$ are not closed: the edge $0\to3$
leaves~$V_1$, and the edges $3\to2$, $3\to5$ leave~$V_2$ into~$V_3$.
The component $V_3$ is the unique closed (absorbing) strongly connected
component of~$G$.
\end{proposition}

\begin{proof}
Immediate from inspection of the edge set~$E$: every edge leaving a
vertex of~$V_3$ lands in~$V_3$ (closure), and every vertex of~$V_3$
is reachable from every other via the edges listed in
Proposition~\ref{prop:transiciones-mod6} (strong connectivity).
\end{proof}

\begin{remark}
The absorbing property of~$V_3$ is the structural reason that odd
integer orbits are confined to~$G'=G[V_3]$ from level~$2$ onwards,
as will be shown in Section~\ref{sec:paths}.
\end{remark}

\begin{proposition}%
\label{prop:ciclos-simples}
The simple cycles of~$G$ are exactly:
\begin{itemize}
  \item Length~$1$: $\;0\to0$ \quad and \quad $5\to5$.
  \item Length~$2$: $\;2\to1\to2$ \quad and \quad $2\to4\to2$.
  \item Length~$3$: $\;2\to1\to5\to2$ \quad and \quad $2\to4\to5\to2$.
\end{itemize}
In particular, $G$ contains no simple cycle of length~$\ge4$.
\end{proposition}

\begin{proof}
The self-loops follow directly from~$E$. Any simple cycle of length~$\ge2$
lies entirely within~$V_3$, since $V_1$ and $V_2$ are not closed.

Within~$V_3$, vertex~$4$ is the key obstruction: its only outgoing edges
are $4\to2$ and $4\to5$, and neither~$1$ nor~$5$ has an edge to~$4$
(only $2\to4$ enters it). Hence any simple cycle through~$4$ enters
via $2\to4$, and leaves either by $4\to2$ (giving $2\to4\to2$) or by
$4\to5$, whence the only return without repeating a vertex is $5\to2$
(giving $2\to4\to5\to2$); no longer simple cycle can contain~$4$.

A cycle avoiding~$4$ lives in $\{1,2,5\}$, where the edges are
$1\to\{2,5\}$, $2\to1$, $5\to\{2,5\}$. Direct enumeration yields
exactly $2\to1\to2$ and $2\to1\to5\to2$. Together with the two cycles
through~$4$, the six listed cycles are the only simple cycles of~$G$.
\end{proof}

\begin{proposition}%
\label{prop:spectral}
Let $A_{G'}$ denote the adjacency matrix of the effective subgraph~$G'$,
and let $A_{\mathcal{H}_4}$ denote the adjacency matrix of the avoiding
subgraph $\mathcal{H}_4:=G'[\{1,2,5\}]$ (obtained by removing vertex~$4$
from~$G'$). Denoting by $\rho(\cdot)$ the spectral radius (largest
absolute value of the eigenvalues), we have
\[
\rho(G')=2,\qquad \rho(\mathcal{H}_4)=\varphi=\tfrac{1+\sqrt{5}}{2},
\]
and the characteristic polynomials factor completely as
\[
\det(\lambda I - A_{G'}) = \lambda^2(\lambda-2)(\lambda+1),
\qquad
\det(\lambda I - A_{\mathcal{H}_4}) = \lambda(\lambda^2-\lambda-1).
\]
The ratio $\varphi/2<1$ is the exact exponential decay rate of the
Fibonacci density, as shown in Section~\ref{sec:paths}.
\end{proposition}

\begin{proof}
In vertex order $(1,2,4,5)$, the edges of~$G'$ are
$1\to\{2,5\}$, $2\to\{1,4\}$, $4\to\{2,5\}$, $5\to\{2,5\}$
(Proposition~\ref{prop:transiciones-mod6}), giving adjacency matrix
\[
A_{G'}=
\begin{pmatrix}
0 & 1 & 0 & 1\\
1 & 0 & 1 & 0\\
0 & 1 & 0 & 1\\
0 & 1 & 0 & 1
\end{pmatrix}.
\]
A direct computation gives $\det(\lambda I - A_{G'})=\lambda^2(\lambda-2)(\lambda+1)$,
with eigenvalues $0,0,2,-1$, hence $\rho(G')=2$.

In vertex order $(1,2,5)$, removing vertex~$4$ gives
\[
A_{\mathcal{H}_4}=
\begin{pmatrix}
0 & 1 & 1\\
1 & 0 & 0\\
0 & 1 & 1
\end{pmatrix}.
\]
A direct computation gives $\det(\lambda I - A_{\mathcal{H}_4})=\lambda(\lambda^2-\lambda-1)$,
whose positive root is $\varphi=(1+\sqrt{5})/2$, so $\rho(\mathcal{H}_4)=\varphi$.
\end{proof}

These three propositions supply all graph-theoretic input needed for
Section~\ref{sec:paths}.

%% file: sec4_block_structure.tex
\section{Modular path structure and the Fibonacci theorem}
\label{sec:paths}

The following theorem describes the exact action of the displacement~$+2^m$
on every level of the modular encoding simultaneously.

\begin{theorem}[Modular path theorem]%
\label{thm:camino-modular}
Let $m\ge1$, $n\in\mathbb{Z}$, and $1\le i\le m$.
Set $d:=(m+i)\bmod2$. Then:
\begin{enumerate}
  \item[\emph{(i)}] \emph{Displacement $+2^m$.}
    \[
    h_i(n+2^m)\equiv
    \begin{cases}
      h_i(n)\pmod6 & \text{if }2^{i-1}\nmid n,\\[4pt]
      h_i(n)+2\pmod6 & \text{if }2^{i-1}\mid n\text{ and }d=0,\\[4pt]
      h_i(n)+4\pmod6 & \text{if }2^{i-1}\mid n\text{ and }d=1.
    \end{cases}
    \]
  \item[\emph{(ii)}] \emph{Displacement $+2^{m+1}$.}
    \[
    h_i(n+2^{m+1})\equiv
    \begin{cases}
      h_i(n)\pmod6 & \text{if }2^{i-1}\nmid n,\\[4pt]
      h_i(n)+4\pmod6 & \text{if }2^{i-1}\mid n\text{ and }d=0,\\[4pt]
      h_i(n)+2\pmod6 & \text{if }2^{i-1}\mid n\text{ and }d=1.
    \end{cases}
    \]
\end{enumerate}
In particular, $h_i(n+2^m)=h_i(n)$ if and only if $2^{i-1}\nmid n$.
\end{theorem}

\begin{proof}
For $i=1$ the claim is immediate: $2^0\mid n$ always, and
$h_1(n+2^m)\equiv h_1(n)+2^m\pmod6$ with $2^m\equiv2$ ($m$ odd) or $4$
($m$ even), matching~(i); the displacement $2^{m+1}$ gives~(ii) likewise.
For $i\ge2$, apply the translation lemma (Lemma~\ref{lem:traslacion}) at level~$i-1$ with displacement $2^m$:
\[
T^{i-1}(n+2^m)-T^{i-1}(n)=2^{m-i+1}\cdot3^{s_{i-1}(n)},
\]
where $s_{i-1}(n)$ is the Hamming weight of the binary encoding $C^{i-1}(n)$.
Note that $s_{i-1}(n)=0$ if and only if $c_1=\cdots=c_{i-1}=0$, i.e.\
$n,T(n),\ldots,T^{i-2}(n)$ are all even, which holds if and only if
$2^{i-1}\mid n$. Hence $s_{i-1}(n)\ge1$ if and only if $2^{i-1}\nmid n$.

\smallskip\noindent
\emph{Case $s_{i-1}(n)\ge1$, i.e.\ $2^{i-1}\nmid n$.}
Since $3^k\equiv3\pmod6$ for all $k\ge1$, we have
$2^{m-i+1}\cdot3\equiv0\pmod6$ for $i\le m$ (as $m-i+1\ge1$).
Hence $h_i(n+2^m)\equiv h_i(n)\pmod6$.

\smallskip\noindent
\emph{Case $s_{i-1}(n)=0$, i.e.\ $2^{i-1}\mid n$.}
Then $3^0=1$ and the difference equals $2^{m-i+1}$.
Since $(m-i+1)+(m+i)=2m+1$ is odd, the parities of $m-i+1$ and $m+i$ are
always opposite. Therefore:
\[
2^{m-i+1}\bmod6=
\begin{cases}
2 & \text{if }m-i+1\equiv1\pmod2,\quad\text{i.e.\ }d\equiv0\pmod2,\\
4 & \text{if }m-i+1\equiv0\pmod2,\quad\text{i.e.\ }d\equiv1\pmod2,
\end{cases}
\]
giving case~(i). For case~(ii), apply the translation lemma with
displacement $2^{m+1}$: the difference is
$2^{m-i+2}\cdot3^{s_{i-1}(n)}$. When $2^{i-1}\nmid n$ the argument
is unchanged (the difference is $\equiv0\pmod6$). When $2^{i-1}\mid n$
the difference is $2^{m-i+2}=2\cdot2^{m-i+1}$. Since
$2^{m-i+1}\bmod6\in\{2,4\}$, multiplying by~$2$ gives
$2\cdot2\equiv4$ and $2\cdot4\equiv2\pmod6$, swapping the two values.
Hence the shifts $+2$ and $+4$ interchange, giving case~(ii).
\end{proof}

\begin{corollary}%
\label{cor:ciclos}
Let $n$ be an odd integer and $m\ge1$. If $T^m(n)=n$, then
$T^m(n+k\cdot2^m)\ne n+k\cdot2^m$ for all $k\not\equiv0\pmod6$.
\end{corollary}

\begin{proof}
Since $n$ is odd, $s=s(C^m(n))\ge1$ (as $c_1=1$). By
Lemma~\ref{lem:traslacion},
\[
T^m(n+k\cdot2^m)=T^m(n)+k\cdot3^s,
\]
so $h_{m+1}(n+k\cdot2^m)\equiv h_{m+1}(n)+3k\pmod6$ (using $3^s\equiv3\pmod6$
for $s\ge1$).

By Theorem~\ref{thm:camino-modular}(i) with $i=1$ (so $d=(m+1)\bmod2$),
each displacement $+2^m$ shifts $h_1$ by $\delta$ where $\delta=2$ if $d=0$
and $\delta=4$ if $d=1$. Since $n+j\cdot2^m$ is odd for every $j$, the
same shift applies at each step, giving
\[
h_1(n+k\cdot2^m)\equiv h_1(n)+k\delta\pmod6.
\]

The cycle condition $T^m(n)=n$ gives $h_{m+1}(n)=h_1(n)$. For $n+k\cdot2^m$
to be a cycle of length~$m$, we need $h_{m+1}(n+k\cdot2^m)=h_1(n+k\cdot2^m)$,
i.e.\ $k\cdot3\equiv k\cdot\delta\pmod6$, i.e.\ $k(\delta-3)\equiv0\pmod6$.
Since $\delta\in\{2,4\}$, we have $\delta-3\in\{-1,1\}$, and
$\gcd(\delta-3,6)=1$, so this holds if and only if $k\equiv0\pmod6$.
Thus for all $k\not\equiv0\pmod6$, $n+k\cdot2^m$ is not a cycle of length~$m$.
\end{proof}

\begin{remark}
Among any six consecutive elements $n,\,n+2^m,\ldots,n+5\cdot2^m$, at most
one can be a cycle of length~$m$; cycles of length~$m$ are therefore
separated by at least $6\cdot2^m$ in arithmetic progressions of common
difference~$2^m$. This period $6\cdot2^m$ coincides with that of the
modular bijection~$\Psi_m$ (Theorem~\ref{thm:biyeccion-modular}),
reflecting the same underlying periodicity of the graph~$G$.
\end{remark}

\begin{corollary}%
\label{cor:rigidez-2adica}
Let $m\ge1$, $n\in\mathbb{Z}$, and let $v_2(n)=r$ denote the $2$-adic
valuation of~$n$ (i.e.\ the largest integer $r\ge0$ such that $2^r\mid n$). Under the displacement $+2^m$, the levels $i=r+2,\ldots,m$ of
$H^m(n)=(h_1,\ldots,h_m)$ (Definition~\ref{def:codificacion-mod-6}) are preserved:
\[
h_i(n+2^m)=h_i(n)\qquad\text{for all }i\ge r+2.
\]
In particular, if $n$ is odd ($r=0$), the suffix $(h_2,\ldots,h_m)$ is
constant along the entire arithmetic progression $n+k\cdot2^m$,
$k\in\mathbb{Z}$. This rigidity is the key ingredient in the proof of
Theorem~\ref{thm:fibonacci}: the bijection $\phi:n\mapsto(h_2,\ldots,h_m)$
used there is well-defined precisely because this suffix is an intrinsic
invariant of the odd integer~$n$, independent of~$k$.
\end{corollary}

\begin{proof}
By Theorem~\ref{thm:camino-modular}(i), the level~$i$ changes if and only if
$2^{i-1}\mid n$, which holds if and only if $i-1\le r$, i.e.\ $i\le r+1$.
For levels $i\ge r+2$ we have $2^{i-1}\nmid n$, so $h_i(n+2^m)=h_i(n)$.
For $n$ odd, $r=0$, so only level $i=1$ changes; applying the shift $k$
times and noting that the argument holds for any translate, the suffix
$(h_2,\ldots,h_m)$ is constant along the entire progression $n+k\cdot2^m$.
\end{proof}

\begin{corollary}%
\label{cor:subgrafo-efectivo}
Let $n$ be an odd integer and $m\ge2$. Then $h_i(n)\in\{1,2,4,5\}$ for all
$i\ge2$. Equivalently, the vertices~$0$ and~$3$ of~$G$ are never visited by
the modular path of an odd integer at any level $i\ge2$.

Consequently, the modular dynamics of odd integers from level~$2$ onwards is
confined to the subgraph $G':=G[\{1,2,4,5\}]$, which is the unique absorbing
strongly connected component of~$G$ (Proposition~\ref{prop:scc}).
\end{corollary}

\begin{proof}
For odd $n=2k+1$, the map~$T$ gives
\[
T(n)=\frac{3(2k+1)+1}{2}=3k+2,
\]
so $h_2(n)\equiv T(n)\equiv 3k+2\pmod6$, which equals~$2$ if $k$ is even
and~$5$ if $k$ is odd. Hence $h_2(n)\in\{2,5\}$ for every odd~$n$.

We verify that the subgraph induced by $\{1,2,4,5\}$ is closed by inspecting
the transitions of~$G$: the edges leaving each of these vertices are
$1\to2$, $1\to5$, $2\to1$, $2\to4$, $4\to2$, $4\to5$, $5\to2$, $5\to5$,
all of which land in $\{1,2,4,5\}$. We now prove $h_i(n)\in\{1,2,4,5\}$
for all $i\ge2$ by induction on~$i$. The base case $i=2$ gives
$h_2(n)\in\{2,5\}\subset\{1,2,4,5\}$ by the above. For the inductive step,
if $h_i(n)\in\{1,2,4,5\}$, then $h_{i+1}(n)\in\{1,2,4,5\}$ by closure.
\end{proof}

\begin{definition}[Avoiding subgraphs]\label{def:avoiding-subgraphs}
For each vertex $v\in\{1,2,4,5\}$, let $\mathcal{H}_v$ denote the subgraph
of~$G'$ obtained by removing vertex~$v$ and all edges incident to it.
Explicitly:
\begin{itemize}
  \item $\mathcal{H}_4 = G'[\{1,2,5\}]$ with edges
    $1\to2$,\;$1\to5$,\;$2\to1$,\;$5\to2$,\;$5\to5$.
  \item $\mathcal{H}_1 = G'[\{2,4,5\}]$ with edges
    $2\to4$,\;$4\to2$,\;$4\to5$,\;$5\to2$,\;$5\to5$.
  \item $\mathcal{H}_5 = G'[\{1,2,4\}]$ with edges
    $1\to2$,\;$2\to1$,\;$2\to4$,\;$4\to2$.
  \item $\mathcal{H}_2 = G'[\{1,4,5\}]$ with edges
    $1\to5$,\;$4\to5$,\;$5\to5$.
\end{itemize}
\end{definition}

The spectral gap $\rho(\mathcal{H}_4)/\rho(G')=\varphi/2$ is what makes the
Fibonacci count exact rather than approximate. The number of paths of
length~$m$ in the full graph~$G$ avoiding vertex~$4$ \emph{nearly} satisfies
the Fibonacci recurrence, failing by exactly~$1$ at each step due to the
self-loop $0\to0$, which contributes paths starting at vertex~$0$ and hence
outside of~$G'$. Restricting to odd integers, whose orbits enter~$G'$ at
level~$2$ (Corollary~\ref{cor:subgrafo-efectivo}), eliminates vertex~$0$
entirely, and the recurrence becomes \emph{exactly} Fibonacci.

\begin{theorem}[Fibonacci count of vertex-$4$-avoiding paths]%
\label{thm:fibonacci}
Let $F(m)$ denote the Fibonacci sequence with $F(1)=F(2)=1$.
For each $\ell\ge0$, the number of paths of length~$\ell$ in~$\mathcal{H}_4$
that start in $\{2,5\}$ is exactly $F(\ell+3)$.

Consequently, the number of odd integers $n\in\{1,\ldots,2^m\}$ whose
modular path $H^m(n)$ avoids vertex~$4$ at every level $i\ge2$ is
exactly $F(m+1)$ (see Example~\ref{ex:fibonacci} for the case $m=4$).
\end{theorem}

\begin{proof}
\emph{Step 1: Paths in $\mathcal{H}_4$.}
For each vertex $v\in\{1,2,5\}$ and each $\ell\ge0$, let $P_\ell(v)$ denote
the number of paths of length~$\ell$ in~$\mathcal{H}_4$ starting at~$v$.
The edges of $\mathcal{H}_4$ are $1\to\{2,5\}$, $2\to1$, and $5\to\{2,5\}$,
so for $\ell\ge1$
\begin{align*}
P_\ell(2)&=P_{\ell-1}(1),\\
P_\ell(1)&=P_{\ell-1}(2)+P_{\ell-1}(5),\\
P_\ell(5)&=P_{\ell-1}(2)+P_{\ell-1}(5),
\end{align*}
with $P_0(v)=1$ for all~$v$. In particular $P_\ell(1)=P_\ell(5)$ for all
$\ell\ge0$. Write $f_\ell:=P_\ell(2)$ and $g_\ell:=P_\ell(5)$. Then
$f_\ell=g_{\ell-1}$ and
$g_\ell=f_{\ell-1}+g_{\ell-1}=g_{\ell-2}+g_{\ell-1}$, with $g_0=1$ and
$g_1=2$, giving $g_\ell=F(\ell+2)$ and $f_\ell=F(\ell+1)$.
The total number of paths of length~$\ell$ starting in $\{2,5\}$ is
\[
f_\ell+g_\ell=F(\ell+1)+F(\ell+2)=F(\ell+3).
\]

\emph{Step 2: Odd integers.}
For $m=1$ the only odd integer in $\{1,2\}$ is $n=1$, for which the
condition is vacuous, giving the count $F(2)=1$; assume henceforth
$m\ge2$.
By Corollary~\ref{cor:subgrafo-efectivo}, $h_2(n)\in\{2,5\}$ for every
odd~$n$. The suffix $(h_2,\ldots,h_m)$ consists of $m-1$ vertices connected
by $m-2$ consecutive edges in~$G'$, so it determines a path of
length~$m-2$ in~$G'$ starting from~$\{2,5\}$. By
Corollary~\ref{cor:rigidez-2adica}, this suffix depends only on~$n$
and not on which element of the progression $n+k\cdot2^m$ is chosen,
so the map $\phi:n\mapsto(h_2,\ldots,h_m)$ is well defined on odd
integers in $\{1,\ldots,2^m\}$.

\emph{Injectivity of $\phi$.} If $\phi(n)=\phi(n')$, then
$h_i(n)=h_i(n')$ for $i=2,\ldots,m$. Since $h_i\equiv c_i\pmod2$,
we get $c_i(n)=c_i(n')$ for $i=2,\ldots,m$. Both $n,n'$ are odd,
so $c_1(n)=c_1(n')=1$, giving $C^m(n)=C^m(n')$. By
Proposition~\ref{prop:biy-fractal}(B), $n=n'$ in $\{1,\ldots,2^m\}$.

\emph{Avoiding vertex~$4$.} The integer $n$ avoids vertex~$4$ at
every level $i\ge2$ if and only if $h_i(n)\ne4$ for $i=2,\ldots,m$,
which holds if and only if $(h_2,\ldots,h_m)$ is a path in~$\mathcal{H}_4$.
Note that $i=1$ is excluded from the condition: for odd~$n$,
$h_1(n)\in\{1,3,5\}$ and the condition $h_1(n)\ne4$ holds automatically.
Thus the number of odd integers avoiding vertex~$4$ is at most
$f_{m-2}+g_{m-2}=F(m-1)+F(m)=F(m+1)$ by Step~1.

\emph{Surjectivity.} Every path $(h_2,\ldots,h_m)$ of length~$m-2$
in~$\mathcal{H}_4$ starting from~$\{2,5\}$ is realised. Indeed,
the string $(1,c_2,\ldots,c_m)$ with $c_i=h_i\bmod2$ determines
a unique odd $n\in\{1,\ldots,2^m\}$ by bijectivity~(B) with
$(h_j(n))_{j=2}^m = (h_j)_{j=2}^m$ by construction: once $h_2(n)\in\{2,5\}$
is fixed by $c_2$, the relation $c_i\equiv h_i\pmod2$ determines each
subsequent $h_i$ uniquely, since from every vertex of~$G$ there is exactly
one even and one odd successor (Proposition~\ref{prop:transiciones-mod6}).
Hence $\phi$ is a bijection and the count equals $F(m+1)$.
\end{proof}

\begin{example}[Fibonacci count for $m=4$]\label{ex:fibonacci}
Among the $8$ odd integers in $\{1,\ldots,16\}$, the integers $3$, $5$, $13$
visit vertex~$4$ (at positions $h_4=4$, $h_3=4$, $h_3=4$ respectively),
while the remaining five avoid it:
\[
\begin{array}{c|c}
n & (h_2,h_3,h_4) \\\hline
1  & (2,1,2) \\
7  & (5,5,2) \\
9  & (2,1,5) \\
11 & (5,2,1) \\
15 & (5,5,5)
\end{array}
\]
These five suffixes are exactly the five paths of length~$2$ in~$\mathcal{H}_4$
from~$\{2,5\}$, generated by the Fibonacci branching: $f_2+g_2=F(3)+F(4)=2+3=5=F(5)$.

For $m=3$ the four odd integers in $\{1,\ldots,8\}$ give the suffixes
$(h_2,h_3)=(2,1)$, $(5,2)$, $(2,4)$, $(5,5)$ for $n=1,3,5,7$ respectively.
Only $n=5$ meets vertex~$4$, so the count is $3=F(4)$. Note that the
condition constrains $h_2$ and $h_3$, that is the iterates $T(n)$ and
$T^2(n)$, and not $T^3(n)$.
\end{example}

\begin{remark}[Elementary formulation via parity words]
\label{rem:palabras-paridad}
Theorem~\ref{thm:fibonacci} admits an equivalent statement in which no graph
appears. For odd~$n$, the modular path $H^m(n)$ avoids vertex~$4$ at every
level $i\ge2$ if and only if the binary encoding
$C^m(n)=(c_1,\ldots,c_m)$ contains no factor~$00$; that is, if and only if
no two consecutive terms of $n,T(n),\ldots,T^{m-1}(n)$ are even.

Indeed, the only edge of~$G'$ entering vertex~$4$ is $2\to4$
(Definition~\ref{def:avoiding-subgraphs}), and $h_2(n)\in\{2,5\}$ for
odd~$n$ (Corollary~\ref{cor:subgrafo-efectivo}). If $h_i(n)=4$ for some
$i\ge2$, then $i\ge3$ and $h_{i-1}(n)=2$; both vertices are even, so
$c_{i-1}=c_i=0$. Conversely, if $c_{i-1}=c_i=0$ for some~$i$, then $i\ge3$
because $c_1=1$, and $h_{i-1},h_i\in\{2,4\}$; since the even successor
of~$2$ is~$4$ and the even successor of~$4$ is~$2$, either $h_{i-1}=4$ or
$h_i=4$, and in both cases vertex~$4$ occurs at a level~$\ge2$.

Since $\Phi_m$ restricts to a bijection between the odd integers in
$\{1,\ldots,2^m\}$ and the words $(c_1,\ldots,c_m)$ with $c_1=1$
(Proposition~\ref{prop:biy-fractal}(B)), the count reduces to the number
$u_m$ of binary words of length~$m$ that begin with~$1$ and avoid~$00$.
Writing $v_m$ for those beginning with~$0$, one has $u_m=u_{m-1}+v_{m-1}$
and $v_m=u_{m-1}$, hence $u_m=u_{m-1}+u_{m-2}$ with $u_1=1$ and $u_2=2$,
so $u_m=F(m+1)$.

This route uses only the bijectivity of~$\Phi_m$. It is worth recording both
because it makes the statement verifiable by hand for small~$m$ and because
it locates the exact relation to the parity-word enumeration
of~\cite{Winkler2026}, whose forbidden factor is~$11$ --- the condition
$T^k(n)\equiv3\pmod4$ --- rather than~$00$. The two conditions are
complementary and are counted by the same transfer matrix
$\left(\begin{smallmatrix}1&1\\1&0\end{smallmatrix}\right)$.
What the modular graph contributes beyond the count itself is the content of
Theorem~\ref{thm:indispensabilidad} and
Corollaries~\ref{cor:visits-v2} and~\ref{cor:quantitative}: the spectral
hierarchy $1<\sqrt{2}<\varphi<2$ across the four vertices of~$G'$ and the
resulting constraints on positive cycles. These are statements about residues
modulo~$6$ and have no formulation in terms of parity words alone.
\end{remark}

\begin{corollary}%
\label{cor:fibonacci-densidad}
The proportion of odd integers in $\{1,\ldots,2^m\}$ whose modular path
avoids vertex~$4$ at every level $i\ge2$ is
\[
\frac{F(m+1)}{2^{m-1}}\;\sim\;\frac{2\varphi}{\sqrt{5}}\left(\frac{\varphi}{2}\right)^m,
\]
where $\varphi=(1+\sqrt{5})/2$. Since $\varphi/2<1$, this tends to zero
exponentially as $m\to\infty$.
\end{corollary}

\begin{proof}
By the Binet formula, $F(m+1)=(\varphi^{m+1}-\psi^{m+1})/\sqrt{5}$, where
$\psi=(1-\sqrt{5})/2$, so $F(m+1)\sim\varphi^{m+1}/\sqrt{5}$ as
$m\to\infty$ (since $|\psi|<1$). Therefore
$F(m+1)/2^{m-1}\sim(2\varphi/\sqrt{5})\cdot(\varphi/2)^m\to0$ since
$\varphi/2<1$. Equivalently, by Proposition~\ref{prop:spectral}, the decay
rate $\varphi/2$ is the ratio of the spectral radii
$\rho(\mathcal{H}_4)/\rho(G')=\varphi/2$ of the avoiding subgraph
and the effective subgraph.
\end{proof}

\begin{proposition}%
\label{prop:fibonacci-matrix}
Let $F_n$ denote the Fibonacci sequence with $F_0=0$, $F_1=1$, extended
by $F_{-1}=1$
(so $F_n=F(n+1)$ in the notation of Theorem~\ref{thm:fibonacci}). In vertex
order $(1,2,5)$, the $n$-th power of the adjacency matrix of~$\mathcal{H}_4$
is
\[
A_{\mathcal{H}_4}^{\,n}=
\begin{pmatrix}
F_{n-1} & F_n & F_n\\
F_{n-2} & F_{n-1} & F_{n-1}\\
F_{n-1} & F_n & F_n
\end{pmatrix},
\qquad n\ge1.
\]
In particular:
\begin{enumerate}
\item[\emph{(i)}] Every entry of $A_{\mathcal{H}_4}^{\,n}$ is a Fibonacci
  number and satisfies the Fibonacci recurrence in~$n$.
\item[\emph{(ii)}] $\operatorname{tr}(A_{\mathcal{H}_4}^{\,n})
  =2F_{n-1}+F_n=L_n$, the $n$-th Lucas number.
\item[\emph{(iii)}] The sum of all entries of $A_{\mathcal{H}_4}^{\,n}$
  is~$F_{n+4}$.
\end{enumerate}
The same holds for $A_{\mathcal{H}_1}$ (in vertex order $(2,4,5)$), by
the isomorphism $\mathcal{H}_4\cong\mathcal{H}_1$.
\end{proposition}

\begin{proof}
The adjacency matrix is
\[
A=A_{\mathcal{H}_4}=
\begin{pmatrix}0&1&1\\1&0&0\\0&1&1\end{pmatrix}.
\]
Since rows~$1$ and~$3$ of $A$ are identical, every power $A^n$ inherits
this symmetry: rows~$1$ and~$3$ of $A^n$ are always equal. We prove
the formula by induction on~$n$.

\emph{Base case $n=1$.} Direct inspection gives
\[
A^1=
\begin{pmatrix}0&1&1\\1&0&0\\0&1&1\end{pmatrix}
=
\begin{pmatrix}F_0&F_1&F_1\\F_{-1}&F_0&F_0\\F_0&F_1&F_1\end{pmatrix}
\]
with $F_{-1}=1$, matching the formula.

\emph{Inductive step.} Assume the formula holds at~$n$. Then
\[
A^{n+1}=A\cdot A^n=
\begin{pmatrix}0&1&1\\1&0&0\\0&1&1\end{pmatrix}
\begin{pmatrix}F_{n-1}&F_n&F_n\\F_{n-2}&F_{n-1}&F_{n-1}\\F_{n-1}&F_n&F_n\end{pmatrix}.
\]
Computing each row: rows~$1$ and~$3$ give
\[
(F_{n-2}+F_{n-1},\;F_{n-1}+F_n,\;F_{n-1}+F_n)=(F_n,\;F_{n+1},\;F_{n+1}),
\]
and row~$2$ gives $(F_{n-1},\;F_n,\;F_n)$.
This matches the formula at~$n+1$.

Part~(ii): $\operatorname{tr}(A^n)=2F_{n-1}+F_n=F_{n-1}+F_{n+1}=L_n$
by the standard identity $L_n=F_{n-1}+F_{n+1}$.

Part~(iii): the sum of all nine entries of $A^n$ equals
$4F_n+4F_{n-1}+F_{n-2}=4F_{n+1}+F_{n-2}$.
We verify this equals $F_{n+4}$ by checking the recurrence and two initial
values: at $n=1$, the sum is $4F_2+F_{-1}=4+1=5=F_5$; at $n=2$, it is
$4F_3+F_0=8+0=8=F_6$. Since both sequences satisfy the Fibonacci recurrence
and agree at $n=1,2$, they are equal for all $n\ge1$.
\end{proof}

\begin{proposition}%
\label{prop:avoiding-v1}
The subgraph $\mathcal{H}_1=G'[\{2,4,5\}]$ is isomorphic to
$\mathcal{H}_4=G'[\{1,2,5\}]$ as directed graphs, and their adjacency
matrices share the same characteristic polynomial
$\lambda(\lambda^2-\lambda-1)$. In particular,
$\rho(\mathcal{H}_1)=\varphi$ and, for each $\ell\ge0$, the number of
paths of length~$\ell$ in~$\mathcal{H}_1$ starting in~$\{2,5\}$ is
$F(\ell+3)$.

Consequently, the number of odd integers $n\in\{1,\ldots,2^m\}$ whose
modular path avoids vertex~$1$ at every level $i\ge2$ is exactly $F(m+1)$.
\end{proposition}

\begin{proof}
The edges of $\mathcal{H}_1$ are $2\to4$, $4\to\{2,5\}$, $5\to\{2,5\}$.
The out-degrees are $\deg^+(2)=1$, $\deg^+(4)=2$, $\deg^+(5)=2$.
The correspondence $1\leftrightarrow4$ (with $2$ and $5$ fixed) defines
a graph isomorphism $\mathcal{H}_4\cong\mathcal{H}_1$: in both subgraphs,
vertex~$2$ has exactly one successor (the removed vertex's counterpart)
and the other two vertices each have successors $\{2,5\}$. The adjacency
matrix of $\mathcal{H}_1$ in vertex order $(2,4,5)$ is
\[
A_{\mathcal{H}_1}=
\begin{pmatrix}
0 & 1 & 0\\
1 & 0 & 1\\
1 & 0 & 1
\end{pmatrix},
\]
which is the same as the adjacency matrix of $\mathcal{H}_4$ in vertex
order $(2,1,5)$, confirming
$\det(\lambda I - A_{\mathcal{H}_1}) = \lambda(\lambda^2-\lambda-1)$.
The path-counting argument is identical to the proof of
Theorem~\ref{thm:fibonacci}, replacing vertex~$1$ by vertex~$4$
throughout (the role of the degree-$1$ vertex at~$2$ is now played
by the edge $2\to4$).
\end{proof}

\begin{remark}
The isomorphism $\mathcal{H}_4\cong\mathcal{H}_1$ reflects a structural
property of~$G'$: vertices $1$ and $4$ are the two successors of
vertex~$2$, and they play symmetric roles in the branching structure.
Removing either one reduces $\deg^+(2)$ from~$2$ to~$1$, creating
the bottleneck that produces the Fibonacci recurrence.
\end{remark}

\begin{proposition}%
\label{prop:avoiding-v5}
The subgraph $\mathcal{H}_5=G'[\{1,2,4\}]$ has adjacency matrix
\[
A_{\mathcal{H}_5}=
\begin{pmatrix}
0 & 1 & 0\\
1 & 0 & 1\\
0 & 1 & 0
\end{pmatrix}
\quad\text{(vertex order $(1,2,4)$),}
\]
with characteristic polynomial $\lambda(\lambda^2-2)$ and spectral radius
$\rho(\mathcal{H}_5)=\sqrt{2}$. The simple cycles contained in
$\mathcal{H}_5$ are exactly $2\to1\to2$ and $2\to4\to2$, both of length~$2$
and both passing through vertex~$2$.
\end{proposition}

\begin{proof}
The edges of $\mathcal{H}_5$ are $1\to2$, $2\to1$, $2\to4$, $4\to2$
(the edges $1\to5$, $4\to5$, $5\to2$, $5\to5$ of~$G'$ are removed
with vertex~$5$). The adjacency matrix is read directly. Its
characteristic polynomial is
$\lambda^3-2\lambda=\lambda(\lambda^2-2)$, with eigenvalues
$0,\pm\sqrt{2}$. The simple cycles are found by inspection: $2\to1\to2$
and $2\to4\to2$ are the only closed paths without repeated vertices.
\end{proof}

\begin{corollary}%
\label{cor:descent-forced-v5}
Every path in~$\mathcal{H}_5$ alternates between $\{1,4\}$ and $\{2\}$.
The only simple closed path in~$\mathcal{H}_5$ that avoids vertex~$4$
is the $2$-cycle $1\to2\to1$ (equivalently, $2\to1\to2$, depending on
starting point), which is the \emph{modular projection} of the known
Collatz cycle $1\leftrightarrow2$ (the sequence of residues
$(h_1,h_2,\ldots)$ that the cycle traces in~$G'$). In particular,
any closed walk in~$\mathcal{H}_5$ that avoids vertex~$4$ is a
repetition of this $2$-cycle.
\end{corollary}

\begin{proof}
In~$\mathcal{H}_5$, vertices $1$ and $4$ have out-degree~$1$ (both
pointing to~$2$), and vertex~$2$ has out-degree~$2$ (pointing to
$\{1,4\}$). Every path therefore alternates between $\{1,4\}$ and
$\{2\}$: after visiting~$2$, the path goes to $1$ or $4$; from either,
it returns to~$2$.

A closed path of even length~$2k$ starting from~$2$ makes $k$
choices (at each visit to~$2$, choosing successor $1$ or~$4$). The only
such path
that avoids vertex~$4$ entirely is $2\to1\to2\to1\to\cdots\to2$,
which as a simple closed path is the $2$-cycle $1\to2\to1$
(the modular projection of the Collatz cycle $1\leftrightarrow2$,
traversed starting from vertex~$1$). Every other closed path visits
vertex~$4$ at least once.
\end{proof}

\begin{theorem}[Vertex indispensability in~$G'$]%
\label{thm:indispensabilidad}
Let $v\in\{1,2,4,5\}$ be any vertex of the effective subgraph~$G'$.
\begin{enumerate}
\item[\emph{(i)}]
  The spectral radius of the avoiding subgraph satisfies
  $\rho(\mathcal{H}_v)<\rho(G')=2$. Explicitly:
  \[
  \rho(\mathcal{H}_4)=\rho(\mathcal{H}_1)=\varphi\approx1.618,
  \quad
  \rho(\mathcal{H}_5)=\sqrt{2}\approx1.414,
  \quad
  \rho(\mathcal{H}_2)=1.
  \]
  The inequality $1<\sqrt{2}<\varphi<2$ governs the hierarchy of
  decay rates.
\item[\emph{(ii)}]
  No positive cycle of~$T$ can have its modular projection entirely
  contained in~$\mathcal{H}_2$. Equivalently, every positive cycle
  of~$T$ visits vertex~$2$ in its modular projection.
\item[\emph{(iii)}]
  The proportion of odd integers $n\in\{1,\ldots,2^m\}$ whose modular
  path avoids vertex~$v$ at every level $i\ge2$ tends to zero
  exponentially at rate $\rho(\mathcal{H}_v)/2<1$.
\end{enumerate}
\end{theorem}

\begin{proof}
Recall that $\mathcal{H}_v$ denotes the subgraph of~$G'$ obtained by
removing vertex~$v$ (Definition~\ref{def:avoiding-subgraphs}).

\emph{Part~(i).} The cases $v=4$, $v=1$, $v=5$ follow respectively from
Theorem~\ref{thm:fibonacci},
Proposition~\ref{prop:avoiding-v1},
Proposition~\ref{prop:avoiding-v5} and Corollary~\ref{cor:descent-forced-v5}.
For $v=2$, the subgraph $\mathcal{H}_2=G'[\{1,4,5\}]$ has edges
$1\to5$, $4\to5$, $5\to5$, whose only cycle is the self-loop $5\to5$;
hence $\rho(\mathcal{H}_2)=1$.

\emph{Part~(ii).} Since vertices~$1$ and~$4$ have no incoming edges
in~$\mathcal{H}_2$ (the only edges entering them in~$G'$, namely
$2\to1$ and $2\to4$, are removed with vertex~$2$), any closed walk
in~$\mathcal{H}_2$ must consist entirely of the self-loop $5\to5$.
But $n\equiv5\pmod6$ with $n>0$ gives $T(n)=(3n+1)/2>n$, so $T$
strictly increases along the walk, contradicting periodicity.

\emph{Part~(iii).} By the Perron--Frobenius theorem, the number of paths
of length~$m$ in~$\mathcal{H}_v$ from~$\{2,5\}$ grows as
$\Theta(\rho(\mathcal{H}_v)^m)$, while the total number of paths in~$G'$
from~$\{2,5\}$ grows as $\Theta(2^m)$. The ratio is
$\Theta((\rho(\mathcal{H}_v)/2)^m)\to0$.
\end{proof}

\begin{remark}
Among the four vertices of~$G'$, vertex~$2$ has the largest effect on
the spectral radius: its removal reduces $\rho$ from~$2$ to~$1$, while
removing vertex~$1$ or~$4$ reduces it to~$\varphi$, and removing
vertex~$5$ to~$\sqrt{2}$. The inequality $1<\sqrt{2}<\varphi<2$ ranks
the vertices by their contribution to the path-counting growth rate.
\end{remark}

\begin{corollary}%
\label{cor:visits-v2}
Every cycle of positive integers under~$T$ visits vertex~$2$ in its
modular projection. Equivalently, no positive cycle can have all its
elements congruent to $1$, $4$, or $5$ modulo~$6$.
\end{corollary}

\begin{proof}
If a positive cycle avoids vertex~$2$, its projection lies
in~$\mathcal{H}_2$, contradicting
Theorem~\ref{thm:indispensabilidad}(ii).
\end{proof}

\begin{proposition}%
\label{prop:flow-conservation}
Let $n$ be a positive integer with $T^m(n)=n$, let $s=s(C^m(n))$ denote
the Hamming weight of the binary encoding, and let $n_v$ denote the
number of indices $i\in\{1,\ldots,m\}$ with $h_i(n)=v$, for
$v\in\{1,2,4,5\}$. Then:
\begin{enumerate}
\item[\emph{(i)}]
  $n_2 = n_1 + n_4$.
\item[\emph{(ii)}]
  $2\,n_2 + n_5 = m$, \quad $s = n_1 + n_5$, \quad $m - s = n_1 + 2\,n_4$.
\item[\emph{(iii)}]
  The distribution $(n_1,n_2,n_4,n_5)$ is uniquely determined by
  the triple $(m,s,n_1)$:
  \[
  n_5 = s - n_1,\qquad
  n_4 = \frac{m-s-n_1}{2},\qquad
  n_2 = \frac{m-s+n_1}{2}.
  \]
\end{enumerate}
\end{proposition}

\begin{proof}
Every positive cycle lies entirely in~$G'$: it contains an odd element
(a positive cycle cannot consist only of even terms, which strictly
decrease), whose next iterate lies in $\{2,5\}\subset G'$
(Corollary~\ref{cor:subgrafo-efectivo}), and since $G'$ is closed the whole
cycle stays in $G'$; hence $n_1+n_2+n_4+n_5=m$.

\emph{Part~(i).}
In~$G'$, the only edges arriving at vertex~$1$ are $2\to1$, and the
only edges arriving at vertex~$4$ are $2\to4$
(Proposition~\ref{prop:transiciones-mod6}). In a closed walk, the number
of arrivals at each vertex equals the number of departures, so the
number of departures from~$2$ equals the number of arrivals at~$1$
plus the number of arrivals at~$4$: $n_2 = n_1 + n_4$.

\emph{Part~(ii).}
Substituting $n_2 = n_1 + n_4$ into $n_1 + n_2 + n_4 + n_5 = m$ gives
$2n_1 + 2n_4 + n_5 = m$, i.e.\ $2n_2 + n_5 = m$. Since
$h_i\in\{1,5\}$ if and only if $T^{i-1}(n)$ is odd, $s = n_1 + n_5$.
Then $m - s = (n_1 + n_2 + n_4 + n_5)-(n_1+n_5) = n_2 + n_4
= (n_1+n_4)+n_4 = n_1 + 2n_4$.

\emph{Part~(iii).}
Solving $s = n_1 + n_5$ and $m - s = n_1 + 2n_4$ gives
$n_5 = s - n_1$ and $n_4 = (m-s-n_1)/2$, hence
$n_2 = n_1 + n_4 = (m-s+n_1)/2$.
\end{proof}

\begin{corollary}%
\label{cor:quantitative}
Every positive cycle of length~$m$ satisfies
\[
n_2\;>\;\frac{m(1-\log_32)}{2}\;\approx\;0.1845\,m.
\]
In particular, vertex~$2$ accounts for more than $18\%$ of the steps.
\end{corollary}

\begin{proof}
By Proposition~\ref{prop:flow-conservation}(iii), $2n_2=m-s+n_1\ge m-s$,
since $n_1\ge0$. For any positive cycle, \eqref{eq:ciclos} gives $n=R/(2^m-3^s)$, and $n>0$ requires $2^m>3^s$, i.e.\ $s<m\log_32$,
so $m-s>m(1-\log_32)$ and hence
\[
n_2\;\ge\;\frac{m-s}{2}\;>\;\frac{m(1-\log_32)}{2}\;\approx\;0.1845\,m.\qedhere
\]
\end{proof}

%% file: conclusions_paper1.tex
\section{Conclusions and future work}
\label{sec:conclusions}

This paper develops a framework in which the Collatz dynamics is encoded
simultaneously at two levels --- binary and modular --- via the bijections
$\Phi_m$ and $\Psi_m$. The framework is built on treating the corrective
residue $R(C^m(n))$ as a function of the binary encoding $C^m(n)$ alone.
This allows both bijections to satisfy the same structural
properties --- bijectivity~(B), periodicity~(P) and antisymmetry~(A) ---
at their respective scales. The general modular Collatz graph is
defined in~\cite{KarrasWeger2026}; the bijection~$\Psi_m$ and its
structural properties do not appear there.

The central result is that among the odd integers~$n$ in
$\{1,\ldots,2^m\}$, exactly $F(m+1)$ have no iterate $T^i(n)$ with
$1\le i\le m-1$ in the residue class~$4\pmod6$
(Theorem~\ref{thm:fibonacci}). This is not an approximation:
the count is exact and unconditional, forced by the spectral gap
$\varphi/2<1$ between~$\mathcal{H}_4$ and~$G'$
(Corollary~\ref{cor:fibonacci-densidad}). Fibonacci counts attached to the
Collatz map are not new; see~\cite{AlbertGudmundssonUlfarsson2015}
and~\cite{Winkler2026}. To our knowledge, however, no prior work counts
integers by the avoidance of a residue class modulo~$6$, nor identifies
$\varphi$ as the spectral radius of an explicit subgraph of the modular
graph.

The modular path theorem (Theorem~\ref{thm:camino-modular}) gives complete
information on how the displacement~$+2^m$ acts on the modular encoding of
every integer in $\{1,\ldots,6\cdot2^m\}$: level~$i$ is preserved if and
only if $2^{i-1}\nmid n$, and otherwise shifts by $+2$ or $+4$ modulo~$6$
according to the parity of $m+i$. Equivalently, if $v_2(n)=r$ is the
$2$-adic valuation of~$n$, then only the first $r+1$ levels of the modular
path change; all deeper levels are rigid
(Corollary~\ref{cor:rigidez-2adica}). As a consequence, positive cycles of
length~$m$ within any arithmetic progression of common difference~$2^m$
are separated by multiples of $6\cdot2^m$ (Corollary~\ref{cor:ciclos}).

No vertex of~$G'$ is dispensable: removing any one of them reduces the
spectral radius strictly below~$2$, with hierarchy $1<\sqrt{2}<\varphi<2$
(Theorem~\ref{thm:indispensabilidad}). In particular, every positive cycle
of~$T$ must visit residue class~$2\pmod6$, and a flow conservation identity
forces this class to account for more than $18\%$ of its steps
(Corollary~\ref{cor:quantitative}).

Three questions remain open.

The first concerns effectivity. For an odd integer~$n$, let $m_0(n)$
be the first level at which $h_{m_0}(n)=4$, with $m_0(n)=\infty$ if the
orbit never meets residue~$4\pmod6$ --- as for $n=1$, fixed by
$1\leftrightarrow2$. The Fibonacci theorem gives $m_0(n)>m$ for exactly
$F(m+1)$ odd integers in $\{1,\ldots,2^m\}$. Whether every odd $n>2$ has
$m_0(n)<\infty$ --- equivalently, whether no such orbit stays in
$\mathcal{H}_4$ forever --- is open; an explicit function $f(n)$ with
$m_0(n)\le f(n)$ for all odd $n>2$ would settle it and connect the count
to a quantitative descent condition.

The second question asks whether every positive cycle other than
$1\leftrightarrow2$ must visit all four vertices of~$G'$.
Theorem~\ref{thm:indispensabilidad} establishes this for vertex~$2$.
Proving it for $v\in\{1,4,5\}$ requires controlling the integrality
of $R/(2^m-3^s)$ (see~\eqref{eq:ciclos}) over all closed walks in
each~$\mathcal{H}_v$, a problem that Proposition~\ref{prop:flow-conservation}
reduces to the analysis of a single free parameter~$n_1$.

The third question asks whether the vertex distribution constraints of
Corollary~\ref{cor:quantitative}, combined with~\eqref{eq:ciclos} and
the structure of~$R$, yield lower bounds on the cycle length~$m$.
The spatial separation of Corollary~\ref{cor:ciclos}, the frequency
bound $n_2>0.1845\,m$, and the integrality of~\eqref{eq:ciclos} together
constrain any hypothetical positive cycle from three independent
directions; whether their combination suffices to derive length
bounds remains open.